\documentclass{amsart}
\usepackage{amsfonts}
\usepackage{stmaryrd}
\usepackage{mathrsfs}
\usepackage{graphicx}

\begin{document}

\newtheorem{thm}{Theorem}[section]
\newtheorem{lem}[thm]{Lemma}
\newtheorem{cor}[thm]{Corollary}
\newtheorem{pro}[thm]{Proposition}
\newtheorem{exm}[thm]{Example}

\theoremstyle{definition}
\newtheorem{defn}{Definition}[section]

\theoremstyle{remark}
\newtheorem{rmk}{Remark}[section]

\def\square{\hfill${\vcenter{\vbox{\hrule height.4pt \hbox{\vrule
width.4pt height7pt \kern7pt \vrule width.4pt} \hrule height.4pt}}}$}
\def\T{\mathcal T}

\newenvironment{pf}{{\it Proof:}\quad}{\square \vskip 12pt}

\title{Degenerating slopes with respect to Heegaard distance}


\author{Jiming Ma}

\address{School of Mathematics Science, Fudan University, Shanghai,
China,
  200433}

\email{majiming@fudan.edu.cn}

\author{Ruifeng Qiu}
\address{School of Mathematical Science  \\Dalian University of
     Technology\\
Dalian, China\\116024} \email{qiurf@dlut.edu.cn}

\thanks{The authors were supported in part by NSFC}


\begin{abstract}

Let $M=H_{+}\cup_{S} H_{-}$ be a genus $g$ Heegaard splitting with
Heegaard distance $n\geq \kappa+2$: (1) Let $c_{1}$, $c_{2}$ be two
slopes in the same component of $\partial_{-}H_{-}$, such that the
natural Heegaard splitting $M^{i}=H_{+}\cup_{S} (H_{-}\cup_{c_{i}}
2-handle)$ has distance less than $n$, then the distance of $c_{1}$
and $c_{2}$ in the curve complex of $\partial_{-}H_{-}$ is at most
$3\mathfrak{M}+2$, where $\kappa$ and $\mathfrak{M}$ are constants
due to Masur-Minsky. (2) Let $M^{*}$ be the manifold obtained by
attaching a collection of handlebodies $\mathscr{H}$ to
$\partial_{-} H_{-}$ along a map $f$ from $\partial \mathscr{H}$ to
$\partial_{-} H_{-}$. If $f$ is a sufficiently large power of a
generic pseudo-Anosov map, then the distance of the Heegaard
splitting $M^{*}=H_{+}\cup (H_{-}\cup_{f} \mathscr{H})$ is still
$n$. The proofs rely essentially on Masur-Minsky's theory of curve
complex.
\end{abstract}

\maketitle

\vspace*{0.5cm} {\bf Keywords}: Heegaard distance, Dehn filling,
handle addition.\vspace*{0.5cm}



\section{Introduction}

 A Heegaard splitting of a compact orientable $3$-manifold $M$ is a
 decomposition of it along an orientable embedded closed surface $S$ into
  two compression bodies $H_{+}$ and
 $H_{-}$ \cite{sch}, and as  a conscious extension of
A. Casson and C. Gordon's notion of strong irreducibility, J.
Hempel \cite{h01} defined the Heegaard distance of a Heegaard
splitting in terms of the curve complex $\mathscr{C}(S)$ of $S$,
the Heegaard distance is the minimal distance between two curves
$\alpha_{+}$ and $\alpha_{-}$ in $\mathscr{C}(S)$ which bound
disks in $H_{+}$ and $H_{-}$ respectively. If a manifold has a
Heegaard splitting with distance at least $3$ , then it is
irreducible, atoroidal, $\partial-$irreducible, anannular and it
is not a Seifert manifold due to Kobayashi and Hempel \cite{h01},
so by Perelman's proof of Geometrization conjecture of Thurston,
the manifold is hyperbolic.

It is expected that high distance splitting has rigidity
properties, Namazi \cite{n} showed that if a 3-manifold $M$ has a
Heegaard splitting with large distance, then the manifold has
finite mapping class group, which is also predict by
Geometrization conjecture. For other highly interesting
construction of high distance Heegaard splittings, see \cite{lm09}
and \cite{lm07}, and see also \cite{l} for related topics.

Thurston's  hyperbolic Dehn surgery theorem, see \cite{t88} and
 \cite{pp}, says that for a noncompact finite volume complete hyperbolic
3-manifold, to each cusp, all but finitely many Dehn surgery
resulting in a hyperbolic 3-manifold. It has been generalized to
hyperbolic manifolds with totally geodesic boundaries by
Scharlemann-Wu and  Lackenby, see \cite{lac02} and \cite{sw}.

Let $M=H_{+}\cup_{S} H_{-}$ be a Heegaard splitting of a
$3-$manifold with boundary and the Heegaard distance is $n$, we
assume that $H_{+}$ is a handlebody and $H_{-}$ is a compression
body. Let $c$ be a slope in $\partial_{-}H_{-}$, then adding a
$2-$handle along $c$ we obtain a compression body $H^{c}$ and a
Heegaard splitting of the new manifold $M^{c}=H_{+}\cup_{S} H^{c}$,
then we consider the problem on the degeneration of distance of the
Heegaard splitting, it is obvious that its distance is at most $n$.
If the distance is less $n$, then we say $c$ is a \emph{degenerating
slope}, we want to know, is there any non-degenerating slope?

Inspired by  canonical Dehn surgery and handle addition theory, i.e,
Gordon  \cite{g} and Scharlemann-Wu \cite{sw}, we have the following
theorem about Heegaard distance:

\begin{thm}

There are constants $\kappa$ and $\mathfrak{M}$ depending on $g$,
such that if $M=H_{+}\cup_{S} H_{-}$ is a genus $g$ Heegaard
splitting of a bordered $3-$manifold with Heegaard distance $n\geq
\kappa+2$, and $c_{1}$, $c_{2}$ are two degenerating  slopes in
the same component $F$ of $\partial_{-}H_{-}$,
 then the distance of $c_{1}$ and $c_{2}$ in  $\mathscr{C}(F)$ is bounded above by $3\mathfrak{M}+2$.

\end{thm}

The constants $\kappa$ and $\mathfrak{M}$ in the theorem  are due to
Masur-Minsky.

Since the curve complex has infinite diameter by Kobayashi-Luo or
Masur-Minsky, see Proposition 4.6 of \cite{mm99}, we have:

\begin{cor}

If  $M=H_{+}\cup_{S} H_{-}$ is a genus $g$ Heegaard splitting with
Heegaard distance $n\geq \kappa+2$ of a bordered $3-$manifold, then
there are infinite ways to attaching handlebodies to the boundary of
$M$, so that the resulting Heegaard splitting has distance $n$.
\end{cor}

\begin{rmk} In canonical Dehn surgery theory and handle addition theory,
 e.g. \cite{g} and \cite{sw},
two "degenerating curves" are related by the intersection number.
In the curve complex there is an up bound of distance in terms of
the  intersection number, see \cite{mm99} or \cite{b}, but in
general, there is no lower bound depending only on the
intersection number, so our theorem is formulated by distance in
the curve complex, and hence we can not obtain the bounded
cardinality of degenerating curves in the Dehn surgery case but
the bounded diameter. It is the disadvantage of our theorems, but
we have the following:

 \end{rmk}

 Note that the  curve complex of a torus is the well-known Farey graph, see
 \cite{mm99},
 and there is a one-to-one correspondence between the slopes in a torus to the co-prime pairs of
 integers, so analogous to a theorem in \cite{s}, we have:

\begin{thm}
If $\partial M$ is a torus and  $M=H_{+}\cup_{S} H_{-}$ is a
Heegaard splitting with Heegaard distance $n\geq \kappa+2$, then
there is a cone in the up half plane, such that at most one
co-prime lattice in the cone is a non-degenerating slope.

\end{thm}

We also prove a result using pseudo-Anosov theory, that is, in the
spirit of \cite{lac02}:

\begin{thm}

If $M=H_{+}\cup_{S} H_{-}$ is a genus $g$ Heegaard splitting with
Heegaard distance $n \geq \kappa+2$, where $H_{+}$ is a handlebody
and $H_{-}$ is a compression body with
$\partial_{-}H_{-}=F_{1}\sqcup F_{2}\ldots \sqcup F_{l}$, $H_{j}$
is a handlebody which has the same genus as $F_{j}$, fixed a
homeomorphism $\tau_{j}:F_{j}\rightarrow
\partial H_{j}$, for generic pseudo-Anosov map $f_{j}: \partial H_{j}\rightarrow
\partial H_{j}$, there is $n_{j}$  such that if  $H^{*}=H_{-}\cup
_{f_{1}^{m_{1}}\tau_{1}}H_{1}\cup_{f_{2}^{m_{2}}\tau_{2}}H_{2}\ldots
\cup_{f_{l}^{m_{l}}\tau_{l}}H_{l}$ with $m_{j}\geq n_{j}$,
 then the resulting Heegaard splitting
$M^{*}=H_{+}\cup_{S} H^{*}$ has Heegaard distance $n$.
\end{thm}

There are examples that after Dehn filling the Heegaard distance
will degenerate drastically, for example,  Minsky, Moriah and
Schleimer \cite{mms} showed that there are arbitrarily high
Heegaard distance knots, and then do the trivial Dehn surgery, we
get the distance zero Heegaard splittings of $S^3$. Our result
shows that most handle additions and handlebody attachments do not
degenerate the Heegaard distance if the initial distance is large,
this is in the spirit of Thurston's Hyperbolic Dehn surgery and
its generalization, which show that most handle additions and
handlebody attachments of a hyperbolic $3-$manifold result in
hyperbolic $3-$manifolds.

The assumption that the distance is  at least $\kappa+2$ is
essential in our proof, but in some sense, it is also necessary, see
the following example:

\begin{exm}
Let $T$ be a closed torus and  $F$ be a torus with an open disk
removed, then $F \times [0,1]$ is a genus $2$ handlebody, and
attaching a $2-$handle along $\partial F \times 0.5$ we get a
distance $2$ Heegaard spilitting of $T^{2}\times [0,1]$, capping
off any torus boundary of $T^{2}\times [0,1]$ by solid torus we
get a distance zero Heegaard splitting of the solid torus.
\end{exm}

In fact, firstly we ponder that if $M=H_{+}\cup_{S} H_{-}$ is a
unstabilized Heegaard splitting of an irreducible 3-manifold $M$
with boundary, is there any way to capping off the negative
boundary of the compression body so that the resulting Heegaard
splitting of the manifold is unstabilized, and then we consider
the generalized problem on Heegaard distance.

In all of the above theorems, we assume that $H_{+}$ is a handlebody
and $H_{-}$ is a compression body, the proof for the case that both
$H_{+}$ and $H_{-}$ are compression bodies is easy from our proof of
the above theorems. In fact,  we just assume that the distance of
$M=H_{+}\cup_{S} H_{-}$ is $n\geq 3$ in this case, then the results
similar to Theorem 1.1 and Theorem 1.4 can be proved in the same
line.

The paper is organized as follows: We outline some fundamental
results on the curve complex $\mathscr{C}(S)$ which we shall use
in Section 2. In Section 3, based on some theorems by
Masur-Minsky, we prove the Theorem 1.1 and theorem 1.3. In Section
4, build on some well-known facts and Lemma 3.1, we prove Theorem
1.4. In Section 5, we treat the case that  both $H_{+}$ and
$H_{-}$ are compression bodies with non-empty negative boundary.


\section{Preliminaries on curve complex}

For a compact surfaces $F$ of genus at least $1$, Harvey
\cite{h81} defined the curve complex, but we just use the
$1-$skeleton of it, we denote it also by $\mathscr C(F)$, whose
vertices are one-to-one corresponding to the essential
non-peripheral curves in $F$, and two vertices are connected by a
length $1$ arc when they are disjoint in the surface $F$ if $F$ is
not homotopic to a torus or once-punctured torus, and two vertices
are connected by a length $1$ arc if they intersect in one point
in the later cases. (Note that in fact in the torus and
once-punctured torus case the definition is due to \cite{mm99},
which is different from \cite{h81} ).

Let $T$ be a torus, fixed a longitude$-$meridian pair $\lambda-\mu$
of the slopes in  $\mathscr{F}=\mathscr{C}(T)$, it is well-known
then the slopes in  $\mathscr{F}=\mathscr{C}(T)$ are determined by a
pair of co-prime integers, so there is a one-to-one corresponding
between the vertex of $\mathscr{C}(T)$ and $\widehat{Q}=Q\cup
\infty$, i.e, $\lambda$ corresponding to $\infty$ and $\mu$
corresponding to $0$. $\mathscr{C}(T)$ is the well-known Faray
graph, see \cite{mm00}.

Masur and Minsky \cite{mm99} (See also Bowditch \cite{b} and
Hamenstadt \cite{ha07} for other proofs) made the important
progress by showing that the curve complex is hyperbolic in the
sense of Gromov and Cannon.

\begin{thm}
$\mathscr{C}(F)$ is $\delta-$hyperbolic, where  $\delta$ depends
only on the topology  of $F$.
\end{thm}

Let $H$ be a genus $g$ handlebody, and $S$ is its boundary, we
denote by $\mathscr{D}(H)$  the subset of $\mathscr C(S)$ each of
which bounds a disk in $H$.

Recall that a subset $A$ of a metric space $X$ is
\emph{quasi-convex} if there is a constant $\kappa$, such that
$\forall~a,b \in A$, any geodesic $[a,b]$ is in the
$\kappa-$neighborhood of $A$. One key theorem  by Masur and Minsky
 \cite{mm04} is:

\begin{thm}
 The disk set
$\mathscr{D}(H)$ is $\kappa-$quasi-convex in $\mathscr{C}(S)$.
\end{thm}

Let $F$ be a compact surface, $Y$ is an \emph{essential} compact
subsurfaces in $F$, by essential subsurface we  mean that
$\pi_{1}(Y) \rightarrow \pi_{1}(F)$ is injective, and $Y$ has genus
at least $1$.(This definition is from \cite{mm00}, but we just
concern with that $Y$ has genus at least $1$).

Masur and Minsky defined the \emph{subsurface projection}
$\pi_{Y}$ from $\mathscr{C}(F)$ to $\mathscr{P}(\mathscr{C}(Y))$,
the power set of $\mathscr{C}(Y)$: for each vertex $v$ in
$\mathscr{C}(F)$, if  $v \cap Y= \emptyset $, then
$\pi(v)=\emptyset$, otherwise $v \cap Y$ is a curve or a set of
arcs in $Y$, then do surgery with the boundary of $Y$, we get a
set of curves in $Y$, which has diameter at most $2$, and it is
easy to show that  that the subsurface projection map is
$2-$Lipschitz, see Lemma 2.3 of \cite{mm00}.

Another key theorem  due to Masur and Minsky is the following
\cite{mm00}:
\begin{thm}\textbf{(Bounded Geodesic Image)}. Let $Y$ be an essential subsurface of $F$,
 and let $\gamma$ be a geodesic segment, ray, or bi-infinite line in
$\mathscr{C}(F)$, such that $\pi_{Y}(v)\neq \emptyset$ for every
vertex $v$ of $\gamma$. There is a constant $\mathfrak{M}$
depending only on $F$ so that $diam_{\mathscr{C}({Y})}
(\pi(\gamma))\leq \mathfrak{M}$.
\end{thm}

Let $F^{'}$ be a compact surface with one boundary, and $F$ be the
surface obtained by capping off the boundary by a disk, then there
is a natural projection map $P: \mathscr{C}(F^{'})
\longrightarrow\mathscr{C}(F)$ by amalgamating the curves which are
identified up to the disk: if $c_{1}$, $c_{2}$ and the boundary of
$F^{'}$ co-bounded a $3-$punctured sphere, we define
$P(c_{1})=P(c_{2})$ in $\mathscr{C}(F)$. Note that  $P$ is a
distance decreasing map.

\section{The proof of the Theorem 1.1 and Theorem 1.3}

In this section, we prove Theorem 1.1, Corollary 1.2 and Theorem
1.3.

Denoted by $F_{j}, j=1,2,3,\ldots, k $ the components of
$\partial_{-}H_{-}$, then there is a set of disks $\mathcal{E}=\{
E_{1},E_{2},E_{3} \ldots E_{k} \}$ in $H_{-}$ which divides
$H_{-}$ into a handlebody(or a 3-ball)and copies of $ F_{j} \times
I$, we denote the components of $\partial_{+}H_{-}-
\partial \mathcal{E}$ which corresponding to $F_{j}$ by $F^{'}_{j}$,
then each $F^{'}_{j}$ is a genus larger or equal to $1$ surface with
a disk removed, we also assume that $\partial E_{j}=\partial
F^{'}_{j}$, see Figure 1. Fixed a homeomorphism $\eta_{j}$ between
$F^{'}_{j}\cup E_{j}$ and $ F_{j}$, we have the canonical isometry
from $\mathscr{C}(F^{'}_{j}\cup E_{j})$ to $\mathscr{C} (F_{j})$, we
also denote it by $\eta_{j}$, and note that for two curves in
$\mathscr{C} (F^{'}_{j})$, the distance  of them under the
composition $\eta_{j} P_{j}: \mathscr{C} (F^{'}_{j})\rightarrow
\mathscr{C} (F^{'}_{j}\cup E_{j})\rightarrow \mathscr{C} (F_{j})$
does not depend on the homeomorphism $\eta_{j}$, so we simply denote
$\eta_{j} P_{j}$ also by $P_{j}$.

Our crucial observation is the follow:

\begin{lem}
$\pi: \mathscr{D_{+}}\longrightarrow \mathscr{C}({F^{'}_{j}})$ is
bounded with constant $\mathfrak{M}$ depending only on $g(S)$.
\end{lem}

\begin{pf}$\forall ~\alpha \in \mathscr{D_{+}}$, since $ E_{j}$ is an essential
 disk in $H_{-}$, we have $d(\alpha,\partial E_{j})\geq n\geq \kappa+2$.
For any pair $x,y \in \mathscr{D_{+}}$, let $[x,y]$ be a geodesic
in $\mathscr{C}(S)$, then it is in the $\kappa-$neighborhood of
$\mathscr{D_{+}}$ by Theorem 2.2, so any vertex $v$ of $[x,y]$ has
distance at most $\kappa$ from a curve, say $\alpha$, in
$\mathscr{D_{+}}$, we have that $d(v, \partial F^{'}_{j}=\partial
E_{j})\geq 2$. Note that distance at least $2$ means that there
intersect essentially in the Heegaard surface $S$, then we can use
Theorem 2.3, which conclude that
$d_{\mathscr{C}(F_{j}^{'})}(\pi(x),\pi(y))\leq \mathfrak{M} $,
where $\mathfrak{M}$ is a constant depending only on $S$.
 \end{pf}

\begin{center}
\includegraphics[totalheight=4cm]{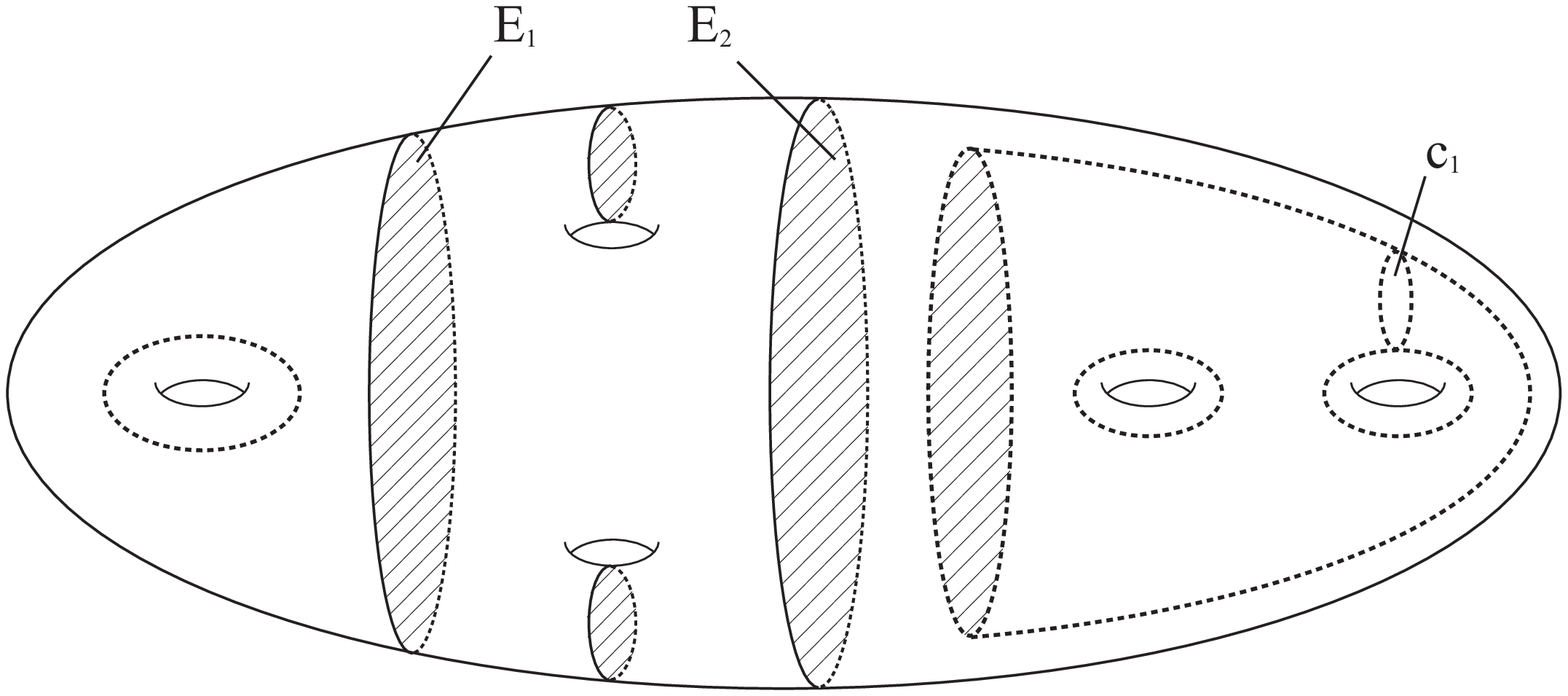}
\begin{center}
Figure 1
\end{center}
\end{center}

\textbf{The proof of Theorem 1.1}:

We assume that two degenerating slopes $c_{1}$ and $c_{2}$ are in
the same boundary component $F=F_{j}$ of $\partial M$,  we also
denote $E_{j}$ by $E$ and $F^{'}_{j}$ by $F^{'}$ for simplicity.
We perform $2-$handle addition along $c_{i}$, and let
$\mathscr{D}_{i}$ be the disk set of
$H_{i}=H_{-}\cup_{c_{i}}2-$handle. We assume that $a_{i}\in
\mathscr{D_{+}}$ and $b_{i}\in \mathscr{D}_{i}$ realize the
Heegaard distance of $M_{i}=H_{+}\cup_{S} H_{i}$, so
$d_{\mathscr{C} (S)} (a_{i},b_{i})\leq n-1$.

Since the Heegaard distance degenerates, we have that $b_{1}\cap
F^{'}\neq \emptyset$, we assume that $b_{1}$ bounds disk $B_{1}$ in
$H_{1}$, and $B_{1}- \mathcal{E}$ is a set of disks, then there is
at least one component of $B_{1}-\mathcal{E} $ which is essential in
$F\times I \cup_{c_{1}}2-handle=H$, we denote the boundary of it by
$c$. In other words,  $c \in \pi(b_{1})$. See Figure 2.

Note that $H$ is a punctured solid torus  or a compression body
according to where the genus of $F$ is $1$ or large. If  $H$ is a
punctured solid torus or
 $g(F)\geq 2$ and $c_{1}$ is seperating in $F$, then there is only one essential disk in
 $H$, its boundary and $c_{1}$ co-bounded an annulus in $F\times
 I$, in this case, we have $d_{\mathscr{C} (F)}(P(c),c_{1})=0$. If $g(F)\geq 2$ and $c_{1}$ is non-seperating in $F$, there is just one non-seperating disk  in
 $H$, and a set of  seperating disk in $H$, the boundary of each separating disk, say $\alpha^{'}$, co-bounded with $\alpha
\subseteq F$ an annulus in $F\times
 I$, and $\alpha$ is disjoint from $c_{1}$, since $\alpha$ is the boundary of a neighborhood of $c_{1} \cup d$, where $d$ is a curve which
 intersects $c_{1}$ with just one point, and we have $d_{\mathscr{C} (F)}(P(c),c_{1})\leq1$.

Let $c^{'}_{1}$ be the curve in $F^{'}$ which co-bounded an annulus
in $F\times I$ with $c_{1}$, note that $d_{\mathscr{C} (S)}
(a_{1},c^{'}_{1})\geq n-1$, suppose otherwise, since $d_{\mathscr{C}
(S)} (\mathscr{D_{-}},c^{'}_{1})\leq 1$, and $d_{\mathscr{C} (S)}
(a_{1},\mathscr{D_{-}})\leq n-1$, a contradiction to the assumption
the initial Heegaard splitting has distance $n$. Let $[a_{1},b_{1}]$
be a geodesic in $\mathscr{C} (S)$, we claim that each vertex of
$[a_{1},b_{1}]$ intersect $F^{'}$ essentially: otherwise, suppose
that $v$ is a vertex which is disjoint from $F^{'}$, then
$d(v,c^{'}_{1})=1$, and we must have that
$d(v,b_{1})=d(v,c^{'}_{1})$ since $a_{1}$ and $b_{1}$ realize the
distance of $M^{1}=H_{+}\cup (H_{-}\cup_{c_{1}}2-handle)$. Then we
have $d(a_{1},v)\leq n-1-1=n-2$, and $v\cap \partial E=\emptyset$,
so $d(a_{1},\partial E)\leq n-1$, which is a contradiction to the
assumption that the initial Heegaard splitting has distance $n$.

\begin{center}
\includegraphics[totalheight=6cm]{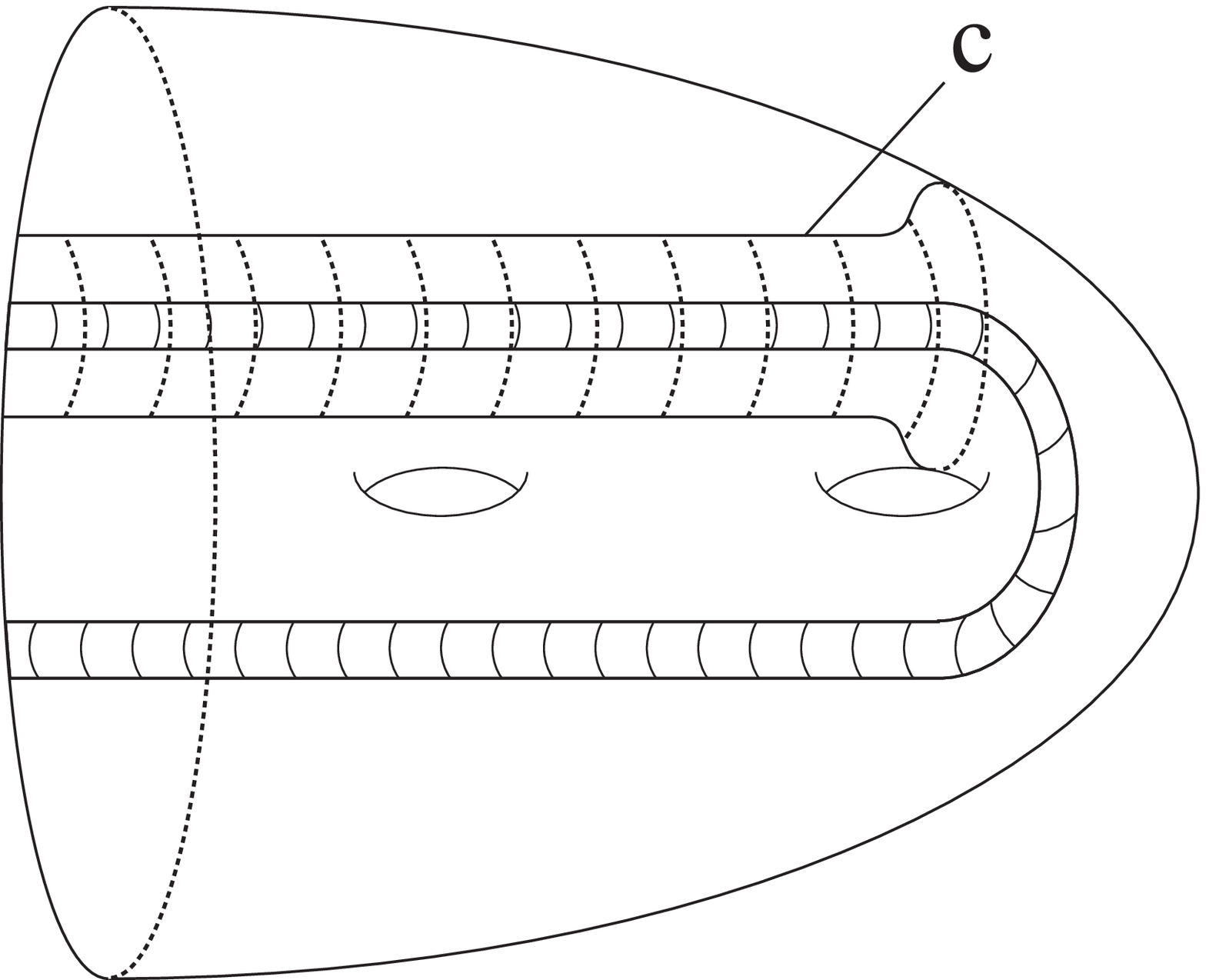}
\begin{center}
Figure 2
\end{center}
\end{center}

From the above claim, and Theorem 2.3, we have $d_{\mathscr{C}
(F^{'})} (\pi(a_{1}),\pi(b_{1}))\leq \mathfrak{M}$. We also claim
that for each vertex $v$ in $[a_{1},b_{1}]$, $P \pi(v)$ is not empty
in $\mathscr{C}(F)$: this is due to the assumption that $F^{'}$ has
just one boundary.


Then since the natural projection map $P$ is distance decreasing,
and each of the projection is not empty, we have $d_{\mathscr{C}
(F)} (P \pi(a_{1}),P \pi (b_{1}))\leq \mathfrak{M}$. Togather with
$d_{\mathscr{C} (F)}(c_{1},P \pi (b_{1}))\leq 1$, we have
$d_{\mathscr{C} (F)}(c_{1},\pi (a_{1}))\leq \mathfrak{M}+1$.
Similarly, we also have $d_{\mathscr{C} (F)} (P \pi(a_{2}),P \pi
(b_{2}))\leq \mathfrak{M}$.

 By Lemma 3.1, $d(P \pi(a_{1}),
P \pi(a_{2}) )\leq \mathfrak{M}$. Then we have $d_{\mathscr{F}
(S)}(c_{1}, c_{2})\leq 3\mathfrak{M}+2$. ~~~~~~~~~~~~~$\Box$

\textbf{The proof of corollary 1.2}:
 This is the  easy corollary of Theorem 1.4, but it also can be
 obtained from Theorem 1.1.

First if some of $F_{j}$ is a genus at least $2$ surface, then
$\mathscr{C} (F_{j})$ is a  diameter infinite graph and each
separating curve $c$ is distance $1$ with a non-separating curve.
So there are infinitely many ways to choose the separating curve
$c_{j}$ and then do $2-$handle addition  along $c_{j}$ we obtain a
distance $n$ Heegaard splitting of a manifold $M^{*}$ with
$\partial M^{*}$ is a set of tori, then for each torus in
$\partial M^{*}$, the curve complex is the Farey graph, which is
also a diameter infinite graph, we have infinite ways to perform
Dehn filling and obtain distance $n$ Heegaard splitting.
$2-$handle additions and Dehn fillings succeeded in the same
boundary is the process of  handlebody attachment, so we perform
handlebody attachments to the manifold with the Heegaard distances
do not degenerate.

Since the distance of $M=H_{+}\cup_{S} H_{-}$ is $n\geq
\kappa+2\geq 3$, so $M$ is hyperbolic with totally geodesic
boundary or with toroidal cusps by Hempel's theorem on Heegaard
distance and Thurston's Hyperbolicity theorem on Haken manifolds,
and $M^{*}$ is a hyperbolic 3-manifolds with toroidal cusps. By
Thurston's Dehn surgery theorem, all but finitely many Dehn
surgery on $M^{*}$ resulting hyperbolic $3-$manifolds with volume
converge to the volume of $M^{*}$. So the manifolds construct
above have infinitely many different volumes, and the handlebody
attachments are different.

 so the non-degenerating slopes
are different even up to homeomorphism of $M$. This means that the
handlebody attachments about are actrally infinitely
many.~~~~$\Box$

 \textbf{The proof of Theorem 1.3}:

Since vertexes  in $\mathscr{F}=\mathscr{C}(T)$ are determined by
a pair of co-prime integers, the pair of co-prime integers are
subset of lattices in $\mathbb{C}$, the plane, but the co-prime
pairs $(a,b)$ and $(-a,-b)$ correspond to the same slope in
$\mathscr{F}=\mathscr{C}(T)$, so we just take lattices in the up
half space $\mathbb{H}$.  We also denote by $b/a$ the lattice
$(a,b)$ in $\mathbb{H}$. Let $L_{\alpha}=\{(x,y)|y=\alpha x\}$ be
a ray in $\mathbb{H}$, and for $\alpha >\theta> 0$, the
$\theta-$neighborhood of $L_{\alpha}$ denoted by
$L_{\alpha}(\theta)$ is the set
$L_{\alpha}(\theta)=\{(x,y)|(\alpha - \theta) x > y > (\alpha
+\theta)x\}$, which is a cone.

We have the following:

\begin{lem}
Let $\mathscr{B}_{n}\subset \mathscr{F}$ be the $n-$neighbourhood
of $1/0$, then there is a cone $L_{\alpha}(\theta)$ such that
 $\mathscr{B}_{n}$ intersect with $L_{\alpha}(\theta)$ by at most
 one point, say $0/1$.
\end{lem}



\begin{pf}We first claim  that for each
$n$, there is a set of cones $L_{\alpha^{n}_{j}}(\theta^{n}_{j})$
which intersect only on $0/1=0$ and three lines
$X_{\pm}=\{(x,y)|x=\pm1\}$ and $Y=\{(x,y)|y=1\}$ with
$\mathscr{B}_{n} \subseteq \cup^{\infty}_{j=1}
L_{\alpha^{n}_{j}}(\theta^{n}_{j}) \cup X_{\pm} \cup Y \cup \{1/0,
0/1\}$.

We prove the claim by induction on $n$. If $n=1$,  note that
$d(1/0,k/1)$ is $1$ in $\mathscr{F}$ since the intersection number
of two slope $b/a$ and $y/x$ is $|bx-ay|$, and so the
$1-$neighborhood of $1/0=\infty$ lies in two vertical lines
$X_{\pm}$ except itself and $0/1$.

The $1-$neighborhood of $0/1$ lies in the horizontal line $Y$
except itself and $1/0$. For a fixed $k\neq 0$, if $d(k/1,b/a)=1$,
then $b=ka\pm1$, so $b/a=k\pm 1/a$ which is in the
$\epsilon_{k}-$neighborhood of $y=kx$ for fixed $\epsilon_{k}> 0$
small enough and $a$ sufficiently large. So there are only
finitely points in the $1-$neighborhood of $k/1$ which is not in
$L_{k}(\epsilon_{k})$, we call them exceptional points, and for
each exceptional point $v/w$, we take a small cone
$L_{v/w}(\epsilon_{v/w})$.

We first choose $\epsilon_{1}$ small enough and then
$\epsilon_{2}$ so  that $L_{1}(\epsilon_{1}) \cap
L_{2}(\epsilon_{2})= \{ 0/1 \}$, and then we choose
$\epsilon_{3},\epsilon_{4},\ldots$. Then we treat the finite
exceptional point $v/w$ for the $1-$neighborhood of $1/1$ one-by
one, choose the $\epsilon_{v/w}$ small enough such the cone
intersects other cones only on $\{ 0/1 \}$, then we treat the
finite exceptional point $p/q$ for the $1-$neighborhood of $2/1$
in the same way, and then the exceptional point for $3/1$...

We reterm these cones by $L_{\alpha^{2}_{j}}(\theta^{2}_{j}),
j=1,2,3\ldots$, and which intersect only on $\{0/1=0\}$ such that
$\mathscr{B}_{2} \subseteq \cup^{\infty}_{j=1}
L_{\alpha^{2}_{j}}(\theta^{2}_{j}) \cup X_{\pm} \cup Y \cup
\{1/0,0/1\}$.

Now assume by induction that the claim is true for $n=k-1$, then
for each $b/a$ which is distance $k-1$ with $1/0$, which is in a
cone constructed above, say in
$L_{\alpha^{k}_{1}}(\theta^{k}_{1})$, choose $\epsilon $ small
enough such that $L_{b/a}(\epsilon) \subset
L_{\alpha^{k}_{1}}(\theta^{k}_{1})$, from the same line above, all
but finitely $1-$neighborhood of $b/a$ lies in
$L_{b/a}(\epsilon)$. For each exceptional point $v/w$ to $b/a$, if
it is also contained by one of the cone constructed above, do
nothing; otherwise, we choose a small cone contains it.

Perform the above process for each point which is distance $k-1$
with $1/0$, we get a set of cones such that $\mathscr{B}_{k}$ are
embraced by these cones together with $X_{\pm} \cup Y \cup
\{1/0,0/1\}$. By induction, the claim follows.

From the claim, we then choose a small cone $L_{\alpha}(\beta)$
which is disjoint from $L_{\alpha^{n}_{j}}(\theta^{n}_{j}),
j=1,2,3\ldots$ but $0/1$, this end the proof of the Lemma.

\end{pf}

 Now from Lemma 3.2 and Theorem 1.1, Theorem 1.3 follows.$\Box$

\section{The proof of the Theorem 1.4}

In this section, we prove Theorem 1.4, first a few facts:

Let $f_{j}: F_{j}\rightarrow F_{j}$ be a pseudo-Anosov map, recall
that for each $f_{j}$, there are two fixed points in
$\mathscr{PMF}(F_{j})$, the projective measured foliation space of
$F_{j}$, say $l_{j+}$ and $l_{j-}$, which are attractor and
repeller respectively, see \cite{flp}. For the fixed handlebody
$H_{j}$, there is the limit set of the mapping class group of
$H_{j}$ acts on $\mathscr{PMF}(F_{j})$, say $\Lambda_{j}$, which
is the closure of the disk set, see \cite{m86}. Since
$\Lambda_{j}$ has measure zero in $\mathscr{PMF}(F_{j})$ by
Kerchhoff \cite{k90}, we say $f_{j}$ is \emph{generic
pseudo-Anosov} if $l_{\pm}\cap \Lambda_{j}=\emptyset $, see
\cite{ns}. Note that generic pseudo-Anosov  $f_{j}$ can not be
extended to a homeomorphism of the handlebody $H_{j}$ and $f_{j}$
acts isometrically on $\mathscr{C}(F_{j})$.


\textbf{The proof of Theorem 1.4}: By Lemma 3.1, the projection of
$\mathscr{D}_{+}$ into $\mathscr{C}({F_{j}})$ has finite diameter.
Let $\mathscr{D}_{j}$ be the set of disks in $H_{j}$, by Theorem
1.1 of \cite{as}, we have that there are two constants $a$ and
$b$, such that in the curve complex of $\mathscr{C}({F_{j}})$, we
have $n/a-b\leq
d_{\mathscr{C}({F_{j}})}(\mathscr{D}_{j},f^{n}_{j}(\mathscr{D}_{j}))\leq
na+b$, so we have that the distance of $\mathscr{D}_{j}$ and
$f^{n}_{j}(\mathscr{D}_{j})$ is large for $n$ sufficiently large,
and then we have
$d_{\mathscr{C}({F_{j}})}(P\pi(\mathscr{D}_{+}),f^{n}_{j}(\mathscr{D}_{j}))$
is large.

 If the Heegaard distance of $M^{*}=H_{+}\cup_{S} H^{*}$ is less
then $n$, we assume that $a$ bounds a disk $\mathscr{D}_{+}$ in
$H_{+}$ and $a^{*}$ bounds a disk $D^{*}$ which realize the Heegaard
distance, then $D^{*}- \mathcal{E}$ is a set of disks, and at least
one, say $D_{j}$, is essential in $H_{j}$. As in the proof of
Theorem 1.1, we have $d(P\pi(\mathscr{D_{+}}), \partial D_{j})$ is
less then $\mathfrak{M}+1$, which is a contradiction to the that
$d(P\pi(\mathscr{D}_{+}),f_{j}(\mathscr{D}_{j}) )$ large. $\Box$


\section{Generalization for $H_{+}$ is a compression body with non-empty negative boundary}

In this section, we generalize the main theorems to the case that
both $H_{+}$ and $H_{-}$ are compression bodies with non-empty
negative boundary. In fact, in this case the assumption on the
initial distance  $n$ can be weaken to $n\geq 3$:

\begin{lem}
If $M=H_{+}\cup_{S} H_{-}$ is a Heegaard splitting with distance
$n\geq 3$, where $H_{+}$ and $H_{-}$ are compression bodies with
non-empty negative boundary. Let $\mathscr{D}_{+}$ be the disk set
of $H_{+}$, then for any component $F$ of $\partial_{-} H_{-}$,
$P\pi(\mathscr{D}_{+})$ has diameter at most $\mathfrak{M}$.

\end{lem}

\begin{pf}Since $\partial_{-} H_{+}$ is non-empty, there is a curve
$c$ in $\mathscr{C}(S)$, which is disjoint from $\mathscr{D}_{+}$.
So $\forall~x,y\in \mathscr{D}_{+}$, we have
$d_{\mathscr{C}(S)}(x,y)=1 ~or ~2$. Note that $x\cap
F^{'}\neq\emptyset$ by the assumption that $n$ is at least $3$,
where $F^{'}$ is the compact surface with one boundary
corresponding to $F$ as in the Section 4, and if $c\cap
F^{'}=\emptyset$, so $d_{\mathscr{C}(S)}(c,\partial F^{'})=1$,
then with $d_{\mathscr{C}(S)}(c,x)=1$ we have
$d_{\mathscr{C}(S)}(x,\partial F^{'})=2$, a contradiction to
$n\geq3$. Then for the length $1$ or $2$ geodesic $[x,y]$, we can
use Theorem 2.3, the lemma follows.
\end{pf}

Now, with Lemma 5.1, similarly to  the proofs in Section 3 and
Section 4, we have:

\begin{thm}

If  $M=H_{+}\cup_{S} H_{-}$ is a genus $g$ Heegaard splitting of a
bordered $3-$manifold with Heegaard distance $n\geq 3$, where
$H_{+}$ and $H_{-}$ are compression bodies with non-empty negative
boundary. If $c_{1}$, $c_{2}$ are two degenerating  slopes in the
same component $F$ of $\partial_{-}H_{-}$,
 then the distance of $c_{1}$ and $c_{2}$ in  $\mathscr{C}(F)$ is bounded above by $3\mathfrak{M}+2$.

\end{thm}

\begin{thm}

If $M=H_{+}\cup_{S} H_{-}$ is a genus $g$ Heegaard splitting with
Heegaard distance $n\geq 3$, where $H_{+}$ and $H_{-}$ are
compression bodies with non-empty negative boundary with
$\partial_{-}H_{-}=F_{1}\sqcup F_{2}\ldots \sqcup F_{l}$, $H_{j}$ is
handlebody which has the same genus as $F_{j}$, fixed a
homeomorphism $\tau_{j}:F_{j}\rightarrow
\partial H_{j}$, for generic pseudo-Anosov map $f_{j}: \partial H_{j}\rightarrow
\partial H_{j}$, there is $n_{j}$  such that if  $H^{*}=H_{-}\cup
_{f_{1}^{m_{1}}\tau_{1}}H_{1}\cup_{f_{2}^{m_{2}}\tau_{2}}H_{2}\ldots
\cup_{f_{l}^{m_{l}}\tau_{l}}H_{l}$ with $m_{j}\geq n_{j}$,
 then the resulting Heegaard splitting
$M^{*}=H_{+}\cup_{S} H^{*}$ has Heegaard distance $n$.
\end{thm}

\end{document}